\documentclass[12pt]{article}
\usepackage{amssymb}
\usepackage{amsmath}
\usepackage{graphicx}
\usepackage{longtable}
\usepackage{psfrag}
\usepackage{mathrsfs}
\parskip \baselineskip
\newtheorem{thm}{Theorem}

\newtheorem{lm}{Lemma}

\newtheorem{prop}{Proposition}

\newtheorem{corol}{Corollary}

\newtheorem{ex}{Example}

\newtheorem{re}{Remark}

\newcommand{\ds}{\displaystyle}

\newcommand{\ol}{\overline}
\newcommand{\ul}{\underline}
\newcommand{\ra}{\rightarrow}
\newcommand{\Ra}{\Rightarrow}
\newcommand{\LRa}{\Leftrightarrow}

\begin{document}
\begin{center}
{\bf Lattices freely generated by posets within a variety. \\
Part I: Four easy varieties} 

Jean Yves Semegni and Marcel Wild
\end{center}

\section{Introduction} \label{section1}

This article constitutes the first of a two part essay triggered by \cite{yves}.  Broadly speaking the intersection of Part I and Part II is the variety of distributive lattices. In Part II this variety is generalized to finitely generated lattice varieties, whereas in the present Part I three other "easy" varieties are taken aboard. These four    varieties are 
\begin{itemize}
\item the variety of all semilattices, 
\item the variety of all lattices, 
\item the variety of all distributive lattices,
\item the variety of all Boolean lattices.  
\end{itemize}
Strictly speaking the first and last are no lattice varieties, yet they fit well.
These four varieties are  easy  in what concerns finding the free object generated by a poset. Some of the semilattice material  will be useful in Part II. All results will be illustrated on a toy poset that accompanies us through both parts. Here comes the section break up. 
\vskip 0.25cm 
An {\it implication} $A\ra B$, where $A, B$ are subsets of some fixed set $P$, is a certain Boolean formula with variables from $P$. In fact, it is the most frequently occuring kind of {\it Horn} formula. If  $\Sigma$ is a family of implications on $P$, then the family $\mathscr{C} (\Sigma)$ of all satisfying truth assignments, or more succinctly of all "$\Sigma$-closed" subsets $X\subseteq P$ is a closure system. Each closure system  $\mathscr{C}$ is of type $\mathscr{C}$ = $\mathscr{C}(\Sigma)$ for a suitable {\it implicational base} $\Sigma$. Section \ref{section2} outlines the $(A,B)$-algorithm \cite{W3} that computes $\mathscr{C}$ from $\Sigma$.
\vskip 0.25cm 
Unless stated otherwise, semilattice means $\vee$-semilattice. According to Theorem 1 in section \ref{section3} each finite presentation $(P,{\cal R})$ of a semilattice $S$ can be rewritten as a family $\Sigma$ of implications on $P$, in which case $S$ turns out to be isomorphic to $\mathscr{C}(\Sigma)\setminus \{\emptyset\}$. As a well known special case one obtains the semilattice $S$ freely generated by a poset $(P,\le)$, which is isomorphic to the semilattice of all nonempty order ideals of $(P,\le)$ under union.
\vskip 0.25cm 
Section \ref{section4} and \ref{section5} are devoted to partial semilattices in general respectively particular. Partial semilattices $(P,\bigvee)$, that is, $\bigvee$ is a certain partial operation from $P^{\omega}$ to $P$, comprise posets $(P,\le)$  as a special case, but are more specific than arbitrary presentations $(P,{\cal R})$ in that the free semilattice $S=F_{\vee}(P,\bigvee)$ generated by $(P,\bigvee)$ cannot collapse $P$. generalizing section \ref{section3}, $S$ can be viewed as a semilattice of certain {\it closed} order ideals of the underlying poset  $(P,\le)$. 
\vskip 0.25cm 
As to the join core  of a lattice $L$ (section \ref{section5}), according to Duquenne \cite{duquenne} it is the unique minimal partial semilattice $(P,\bigvee)$ such that $F_{\vee}(P,\bigvee)\simeq L$ as semilattices. We deemed it worthwhile to reprove this neat fact in a way that relates a bit more to the theory of implicational bases $\Sigma$.
\vskip 0.25cm 
Section \ref{section6} reviews, because it fits well, the lattice $FL(P,\le)$ freely generated by a poset within the variety of all lattices.
 
\vskip 0.25cm 
The structure of the lattice $FD(P,\le)$ freely generated by a poset within the variety ${\cal D}$ of all distributive lattices is also well known, but our proof  in section \ref{section7} seems to be new. 
\vskip 0.25cm 
Albeit as set, the variety ${\cal B}$ of all Boolean lattices is contained in ${\cal D}$, the lattice $FB(P,\le)$ freely generated by $(P,\le)$ within ${\cal B}$ is usually larger that $FD(P,\le)$ because $FB(P,\le)$ additionally needs to be closed under complementation. The structure of $FB(P,\le)$ is determined by the number $t$ of atoms. Finding $t$ in terms of $(P,\le)$ is easy but believed to be new. For our toy poset $(P,\le)$ the cardinalities of $FL(P,\le)$, $FD(P,\le)$, and $FB(P,\le)$ are shown to be $35$, $25$, and $16384$ respectively.

\section{Implications and the $(A,B)$-algorithm} \label{section2}

An {\it implication} on a set $P$ is a pair of subsets $(A,B)$, written as $A \ra B$. Here $A$ and $B$ are the {\it premise} and {\it conclusion} of the implication. Let
\begin{eqnarray}\label{eq1}
 \Sigma & : = & \{A_1 \ra B_1, \cdots, A_n \ra B_n\}
\end{eqnarray}

be a family of implications. A subset $X \subseteq P$ is {\it $\Sigma$-closed} if
\begin{eqnarray}\label{eq2}
(\forall 1 \leq i \leq n) \;\ (A_i \subseteq X \Ra B_i \subseteq X)
\end{eqnarray}
In other words, for each $i$ one must have either $A_i \not\subseteq X$ or $B_i \subseteq X$ (or both). It is easy to see that the family $\mathscr{C}(\Sigma )$ of all $\Sigma$-closed subsets is a {\it closure system}, that is, it contains $P$ and is closed under intersections. Observe that the empty set $\emptyset$ is in $\mathscr{C}(\Sigma)$ if and only if there are no implications of type $\emptyset\ra B_i\; (B_i \neq \emptyset)$. There is no need to deepen the mentioned connection (introduction) to Horn formulae here; that is done in \cite{W3}.

\begin{ex}\label{ex1} {\rm Let} $P = \{a, b, c, d, e, f, g \}$ \; {\rm and} 
\begin{eqnarray*} 
\Sigma  & := & \left\{ \;\{ b\} \ra \{a, e\}, \ \{c\} \ra \{b \}, \ \ \{ d \} \ra \{c, f \}, \right. \\
& & \left. \{f\} \ra \{e\}, \ \{g\} \ra \{b, f\}, \  \{b, f\}  \ra \{g\} \;\right\}
\end{eqnarray*}
\end{ex}

Then $X = \{a, b, e\}$ is $\Sigma$-closed since with $\{b\} \subseteq X$ also $\{a, e\} \subseteq X$. Ditto $\{a, e\}$ is $\Sigma$-closed, but not $\{a, b\}$. Here is the complete closure system:

\newpage

\begin{figure}[!h]
\begin{center}
\psfrag{csig}{$\mathscr{C}(\Sigma)\;\;\;\;\; =$}
\includegraphics[scale=.4]{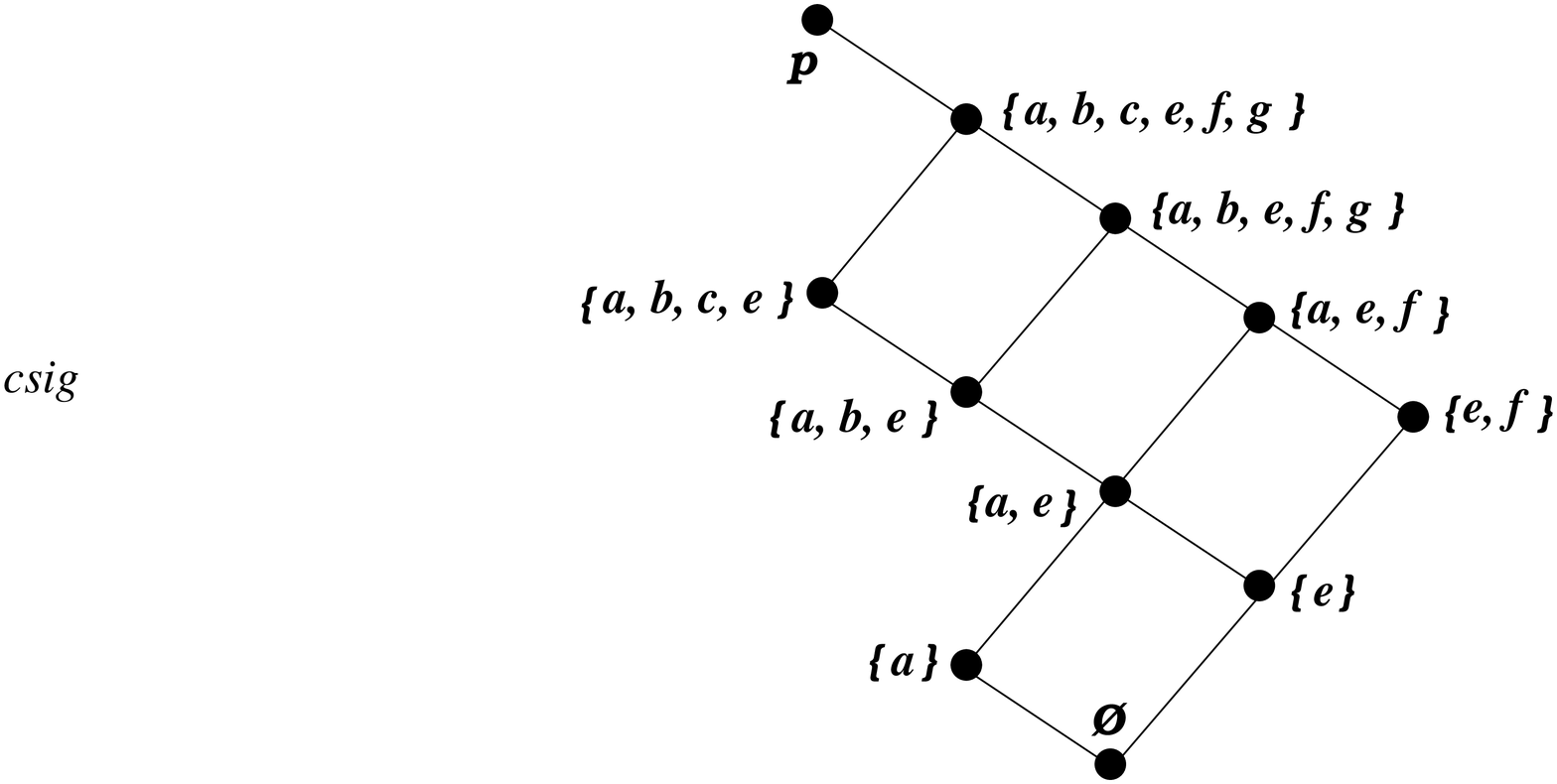}
\caption{}
\label{fig1}
\end{center}
\end{figure}

Let us give a rough sketch of the $(A, B)$-algorithm \cite{W3} which efficiently generates closure systems of the kind $\mathscr{C}(\Sigma)$.

\begin{ex} \label{ex2} \end{ex}
Let $\Sigma$ be as in Example \ref{ex1}. The family $r$ of all sets $X \in 2^P$ that satisfy $\{b\} \ra \{a, e\}$ can compactly be written\footnote{We used $\alpha, \beta$ instead of $a, b$  in \cite{W3} because $a$ and $b$ have another meaning in the present article.} as

\hspace*{2cm}  \begin{tabular}{c|c|c|c|c|c|c|c|c|}
 & $a$ & $b$& $c$ & $d$ & $e$ & $f$ & $g$ \\ \hline
$r=$ & $\beta$  & $\alpha$ & $2$ & $2$ & $\beta$ & $2$ & $2$ \\ \hline \end{tabular}

That is, each $0,1$-incidence vector corresponding to such a set $X$ must be such that if its second component $\alpha$ is $1$, then also $1$ must occur at the two positions labelled $\beta$. The other positions carry a label 2 which indicates that they are free to be independently $0$ or $1$. In order to impose the second implication $\{c\} \ra \{b\}$ we first split $r$ into the disjoint union of $r_1= \{X \in r| \ b \not\in X\}$ and $r_2 = \{ X \in r | \ b\in X\}$. Thus

\hspace*{2cm} \begin{tabular}{c|c|c|c|c|c|c|c|}
           & $a$ & $b$ & $c$ & $d$ & $e$ & $f$ & $g$\\ \hline
$r_1=$ & $\beta$ & $\alpha$ & ${\bf 0}$ & $2$ & $\beta$ & $2$ & $2$ \\ \hline
$r_2=$ & $\beta$ & $\alpha$ & ${\bf 1}$ & $2$ & $\beta$ & $2$ & $2$  \\ \hline \end{tabular}

All $X \in r_1$ satisfy $\{c \} \ra \{b\}$ since $c \not\in X$. But $X \in r_2$ satisfies $\{c\} \ra \{b\}$ if and only if $\alpha =1$ in $r_2$. The latter entails that both $\beta$ in $r_2$ are $1$. Hence the family of all $X \in 2^P$ that satisfy both of $\{b\} \ra \{a,e\}$ and $\{c \} \ra \{b\}$ is

\hspace*{2cm} 
\begin{tabular}{|c|c|c|c|c|c|c|}
$a$ & $b$ & $c$ & $d$ & $e$ & $f$ & $g$\\ \hline
$\beta$ & $\alpha$ & $0$ & $2$ & $\beta$ & $2$ & $2$ \\ \hline
$1$ & $1$ & $1$ & $2$ & $1$ & $2$ & $2$ \\ \hline 
\end{tabular}

Continuing by imposing $\{d\} \ra \{c, f\}$ up to $\{b, f\} \ra \{g\}$ one obtains the table

\hspace*{2cm} 
\begin{tabular}{|c|c|c|c|c|c|c|}
$a$ & $b$ & $c$ & $d$ & $e$ & $f$ & $g$\\ \hline
$\beta$ & $\alpha$ & $0$ & $0$ & $\beta$ & $0$ & $0$\\ \hline
$2$ & $0$ & $0$ & $0$ & $1$ & $1$ & $0$ \\ \hline
$1$ & $1$ & $0$ & $0$ & $1$ & $1$ & $1$ \\ \hline
$1$ & $1$ & $1$ & $0$ & $1$ & $0$ & $0$ \\ \hline
$1$ & $1$ & $1$ & $2$ & $1$ & $1$ & $1$ \\ \hline
\end{tabular}

which encodes $\mathscr{C}(\Sigma)$. For instance, the first row comprises the five sets $\phi, \{a\}, \{e\}, \{a, e\}$  and $\{a, b, e\}$.  Since the rows comprise mutually disjoint set families, one obtains 
\begin{eqnarray*}
|\mathscr{C}(\Sigma)|\;\;\;= \;\;\;5+2+1+1+2 \;\;\;=\;\;\; 11
\end{eqnarray*}
in accordance with Figure \ref{fig1}.

One can show that the wasteful deletion of rows during the $(A, B)$-algorithm can be avoided. Another nice feature is that the time spent is mainly dependent on the number of implications and {\it not} on the actual size of $\mathscr{C}(\Sigma)$. The $(A,B)$-algorithm simplifies to the $(a,B)$-algorithm, and runs particularly well, when all premises are singletons.  

Back from algorithms to theory. As is well known, every closure system $\mathscr{C}$ on $P$ (not necessarily induced by $\Sigma$) comes along with a closure operator $A \mapsto \ol{A} \ (A \subseteq P)$ defined by 
\begin{eqnarray} \label{eq3}
 \ol{A} \quad : =\quad \bigcap \{X \in \mathscr{C} | \ X \supseteq A\}.
\end{eqnarray}
Furthermore $\mathscr{C}$, partially ordered by inclusion, is a $\vee$-semilattice with joins (suprema) given by 
\begin{eqnarray} \label{eq4}
 \quad A \vee B \quad = \quad \ol{A \cup B}.
\end{eqnarray} 

Conversely, every closure operator $A\mapsto \ol{A}$ yields the closure system $\mathscr{C}$ of all {\it closed} sets $A=\ol{A}$. Given any closure system $\mathscr{C}$ on $P$, a family $\Sigma$ of implications with $ \mathscr{C}(\Sigma ) = \mathscr{C}$ is called an {\it implicational base} of $\mathscr{C}$. It is {\it nonredundant} if no proper subset of $\Sigma$ is an implicational base of $\mathscr{C}$. Any two nonredundant implicational bases $\Sigma _1, \Sigma _2$ of $\mathscr{C}$ yield the same set $E(\mathscr{C})$ of {\it essential} elements in that
\begin{eqnarray}\label{eq5}
\{\ol{A}\;|\; (A\ra B)\in \Sigma _1\} &=& \{\ol{C}\;|\; (C\ra D)\in \Sigma _2\}\;\;:=\;\; E(\mathscr{C}).
\end{eqnarray}
An implicational base $\Sigma $ is {\it optimal} if the sum of the cardinalities of all premises and conclusions of implications occuring in $\Sigma$ is minimal. One can show that each optimal implicational base is nonredundant.  A nonredundant implicational base can be computed in quadratic time, but computing an optimal implicational base in $NP$-hard. Interestingly, if $\mathscr{C}$ is modular as a lattice, the task can be achieved in polynomial time.

Along with each closure operator $X\mapsto \ol{X}$ on $P$ comes the {\it quasi-closure} operator $X\mapsto X^{\bullet}$ which is defined by 
\begin{eqnarray} \label{eq6}
 X^{\bullet} &:=& X\cup X^{\circ} \cup X^{\circ \circ} \cup X^{\circ \circ \circ}\cup \cdots, 
\end{eqnarray}
where
\begin{eqnarray}\label{eq7}
Y^{\circ} &:=& Y \cup \bigcup \left\{\ol{Z}\; |\; Z\subseteq Y \; \textrm{and}\; \ol{Z}\ne \ol{Y} \right\}.
\end{eqnarray}

Notice that the iterated sets $X^{\circ \circ \cdots \circ}$ in (\ref{eq6})  eventually become stationary due to the finiteness of $P$. It is clear that $X\mapsto X^{\bullet}$ is a closure operator and that $X^{\bullet}\subseteq \ol{X}$ for all $X\subseteq P$. Call $X\subseteq P$ {\it quasiclosed} if $X^{\bullet}=X$.

\section{Finitely presented semilattices}\label{section3}

Let $P$ be any finite set of ``symbols'' and let ${\cal R}$ be a finite set of semilattice relations with symbols from $P$. The {\it semilattice} $F_\vee (P, {\cal R})$ {\it freely generated by the set $P$ and subject to the relations in} ${\cal R}$ is the up to isomorphism unique \footnote{Recall that unless stated otherwise, "semilattice"  means $\vee$-semilattice. The semilattice $F_\vee (P, {\cal R})$ is a special case of a universal algebra {\it finitely presented by generators and relations.} Such creatures are always unique up to isomorphism.} semilattice $S$ such that

(a) \quad  There is a map $\phi : \ P \rightarrow S$ which {\it satisfies} the relations from $\mathcal{R}$ (obvious definition) \\
\hspace*{1.2cm} and is such that $\phi (P)$ generates $S$.

(b) \quad For each semilattice $T$ and each map $\phi : \ P \ra T$ {\it respecting} the relations from  ${\cal R}$ \\
\hspace*{1.1cm}(obvious definition), there is a $\vee$-homomorphism $\Phi : F_\vee (P, {\cal R}) \ra T$ that extends $\phi$.

Each semilattice relation in ${\cal R}$, say $a\vee b= c\vee d$, can be rewritten as $a\vee b\ge c\vee d$ and $c\vee d\ge a\vee b$. Of course ${\cal R}$ can also feature ``lone'' inequalities $x\ge y$, since this amounts to $x\vee y=x$. So we may suppose that all relations in ${\cal R}$ are of type
\begin{eqnarray} \label{eq8}
\quad a \vee b \vee \cdots \vee c \quad \geq \quad d \vee e \vee \cdots \vee f.
\end{eqnarray}
Let $\Sigma ({\cal R})$ be the family of all corresponding implications
\begin{eqnarray}\label{eq9}
\{a,b, \cdots, c\}\quad  \ra \quad  \{d, e, \cdots, f\}.
\end{eqnarray}
\begin{thm} \label{theorem1} [5, Satz 45]  If all relations in ${\cal R}$ have been adjusted to type (\ref{eq8}), then the semilattice $F_\vee (P, {\cal R})$ is isomorphic to the semilattice $ \mathscr{C} (\Sigma ({\cal R})) \setminus \{\emptyset \}$, with joins as in (4).
\end{thm}

{\bf Proof:} We first show that $S : =  \mathscr{C} (\Sigma ({\cal R})) \setminus \{ \emptyset \}$ satisfies (a). For all $a \in P$ put $\ol{a}: = \ol{\{a\}}$ and define $\phi (a) : = \ol{a}$. That $S$ is generated by $\ol{P} :  = \phi (P)$ can be seen from
$$\ol{\{a, b, \cdots, c\}} \quad = \quad \ol{\ol{a} \cup \ol{b} \cup \cdots \cup \ol{c}} \quad \stackrel{(4)}{=} \quad \ol{a} \vee \ol{b} \vee \cdots \vee \ol{c}.$$
Let $a \vee b \vee \cdots \vee c \geq d\vee e \vee \cdots \vee f$ be a relation from ${\cal R}$. Since $\{a,  \cdots, c\} \ra \{d, \cdots, f\}$ belongs to $\Sigma ({\cal R})$, and $\ol{\{a, \cdots, c\}}$ is $\Sigma ({\cal R})$-closed, we conclude $\ol{\{a, \cdots, c\}} \supseteq \{d, \cdots, f\}$, which yields
$$\ol{a} \vee \cdots \vee \ol{c}\quad =\quad \ol{\{a, \cdots, c\}}\quad \supseteq \quad \ol{\{d, \cdots, f\}}\quad =\quad \ol{d} \vee \cdots \vee \ol{f}$$
This shows (a). As to (b), let $T$ be any semilattice and let $\phi : P \ra T$ respect ${\cal R}$. We claim:
\begin{eqnarray}\label{eq10}
 \textrm{If}\;\;\; g \in \ol{\{a, b, \cdots, c\}}, \;\;\;\textrm{then}  \;\; \;\phi (g) \leq \phi (a) \vee \phi (b) \vee \cdots \vee \phi (c)
\end{eqnarray}
To fix ideas, say $\{a, b\} \ra \{d\}$ and $\{c, d\} \ra \{g\}$ belong to $\Sigma ({\cal R})$, and so $g \in \ol{\{a,b, c\}}$. Since $\phi$ respects the relation $a \vee b \geq d$,  one has 
$\phi (a) \vee \phi (b) \geq \phi (d)$ in $T$. Ditto $\phi (c) \vee \phi (d) \geq \phi (g)$. But this yields $\phi (a) \vee \phi (b) \vee \phi (c) \geq \phi (g)$, which proves (\ref{eq10}). Define $\Phi : S \ra T$ by
$$\Phi (\ol{\{a, b, \cdots, c\}})\quad :=\quad \phi (a) \vee \phi (b) \vee \cdots \vee \phi (c)$$
Because of (\ref{eq10}), $\Phi$ is well defined, i.e. not dependent on the particular generators $a, b, \cdots, c$. From $\ol{\{a, \cdots, c\}} \vee \ol{\{d, \cdots, e\}} = \ol{\{a, \cdots, c, d, \cdots,  e\}}$  readily follows that $\Phi$ is a $\vee$-morphism. \hfill $\blacksquare$

We mention that Theorem \ref{theorem1} generalizes in natural ways to finitely presented commutative semigroups.

\begin{ex} \label{ex3}{\rm By Theorem \ref{theorem1} and Example \ref{ex1}, if $P = \{a, b, \cdots, g \}$, and ${\cal R}$ consists of the  relations $b \geq a \vee e, \ c \geq b, \ d \geq c \vee f, \ f \geq e, \ g = b \vee f$, then}
\end{ex}

\begin{figure}[!h]
\begin{center}
\psfrag{fv}{$F_\vee (P,{\cal R}) \;\;=$}
\includegraphics[scale=.35]{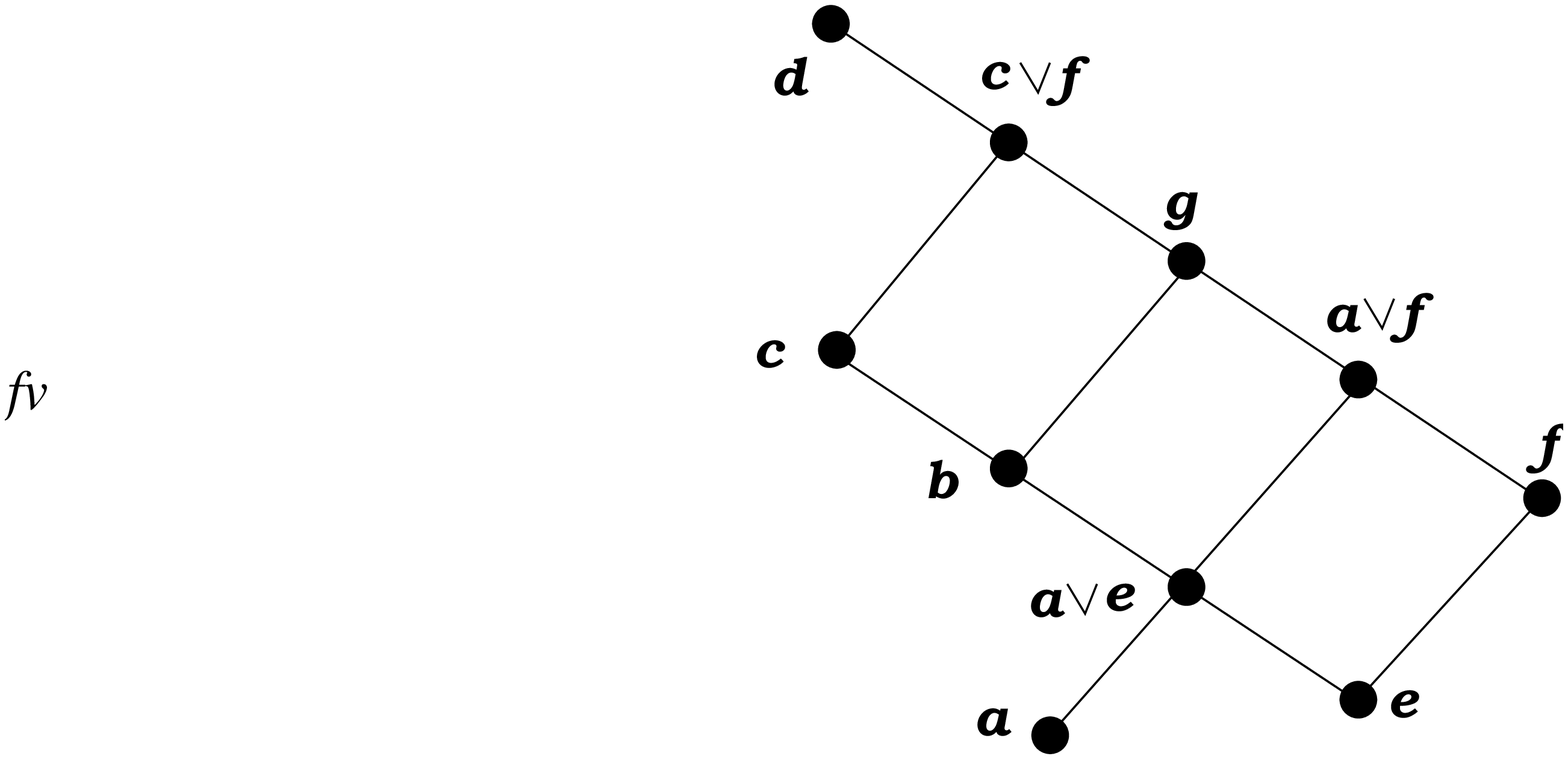}
\caption{}
\label{figure2}
\end{center}
\end{figure}

Relations not present in ${\cal R}$ may well be present in $F_\vee(P,{\cal R})$, such as $e\le c$. But such relations are always {\it deductible}  from the relations in ${\cal R}$. Thus $e\le c$ follows from $e\le a\vee  e\le b\le c$, where the first  inequality follows from the semilattice axioms, and the other two are relations in ${\cal R}$. The relations in ${\cal R}$ may even cause the collapse of $F_\vee(P,{\cal R})$ to a single point, say if ${\cal R} \;\;=\;\;\{a\ge b,\; b\ge c,\; c\ge a\}$.

Let $(P,\le)$ be a finite poset. The semilattice freely generated by $(P,\le)$ is defined as 
$$F_\vee(P,\le)=F_\vee(P,{\cal R})$$
 where ${\cal R}$  consists of all relations $a\ge b$ holding in $P$. Because $\vee$ does not feature in ${\cal R}$, the only possible deductions use the transitivity of $\ge$. But if $a\ge b$ and $b\ge c$ are in ${\cal R}$ then so is $a\ge c$. Hence ${\cal R}$ is deductively closed.
\begin{corol}
The semilattice $F_\vee(P,\le)$ is isomorphic to the $\cup$-semilattice $Id(P,\le)\setminus \{\emptyset\}$. Furthermore, $F_\vee(P,\le)$ is up to isomorphism the unique semilattice $S$ such that:
\begin{itemize}
\item[(i)] There is a generating subset of $S$ which as a poset is isomorphic to $(P,\le)$. 
\item[(ii)]For each semilattice $T$ and each monotone map $\phi : P\ra T$ there is a semilattice morphism $\Phi: S\ra T$ that extends $\phi$. 
\end{itemize}
\end{corol}

{\bf Proof:} By Theorem \ref{theorem1}, $F_\vee(P,\le)$  is isomorphic to $\mathscr{C}\left(\Sigma\left({\cal R}\right)\right)\setminus \{\emptyset\}$ which of course is $Id(P,\le)\setminus \{\emptyset\}$. Properties $(i)$ and $(ii)$ essentially follow from (a), (b) (beginning of sec.\ref{section3}).  The fact that in $(i)$  the generating set of $S$ cannot collapse, but 
is {\it isomorphic} to $P$, is due to the fact that here ${\cal R}$ is deductively closed.  \hfill $\blacksquare$

Mutatis mutandis, all holds when $\vee$ is switched with $\wedge$. Specifically, if $(P^d,\le^d)$ is the dual of $(P,\le)$, then 
\begin{displaymath}
F_\wedge(P,\le) \;\;\cong \;\; F_\vee(P^d,\le^d)\;\; \cong \;\;Id(P^d,\le^d)\setminus \{\emptyset\}\;\;\cong \;\; Fil(P,\le)\setminus\{\emptyset\}, 
\end{displaymath}
where $Fil(P,\le)$ is the family of order filters of $(P,\le)$.

\begin{ex}\label{ex4}
{\rm Along with the poset} $(P,\le)$ {\rm in Fig.3(a)  the freely generated} $\vee${\rm -semilattice} $F_\vee(P,\le)$ {\rm(Fig.3(b))  and the freely generated} $\wedge${\rm -semilattice} $F_\wedge(P,\le)$ {\rm (Fig.3(c))  are shown.} 
\end{ex}

Notice that despite first appearances $F_\wedge(P,\le)$ is {\it not} the dual of $F_\vee(P,\le)$. The above isomorphism $F_\wedge(P,\le) \ra Fil(P,\le)\setminus\{\emptyset\}$ sends (say) $c$ to $c\! \uparrow$ and $f$ to $f\!\uparrow$, but $c\wedge f$ maps to $c\!\uparrow \cup f\!\uparrow$. Thus  the semilattice operation in  $Fil\left(P,\le\right)\setminus\{\emptyset\}$ is $\cup$.

\begin{figure}[!h]
\begin{center}
\psfrag{u}{$(P,\le)\;=\;$}
\psfrag{v}{$F_{\vee}(P,\le)\;=\;$}
\psfrag{w}{$F_{\wedge}(P,\le)\;=\;$}
\includegraphics[scale=.375]{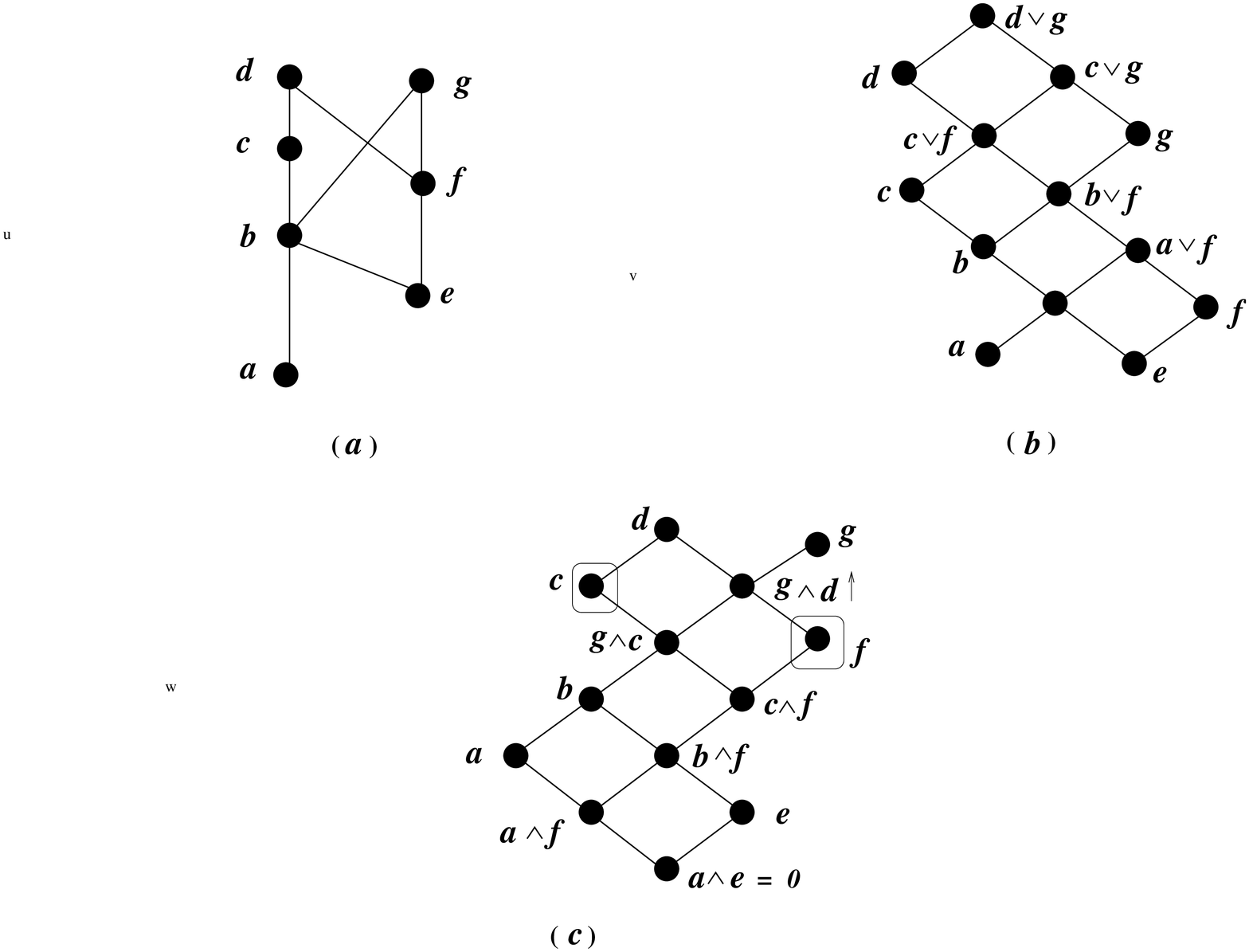}
\end{center}
\caption{}
\label{figure3}
\end{figure}

\section{Partial semilattices}\label{section4}
Every subset $P$ of  a semilattice $S$ determines a partial operation $\bigvee$ from $P^{\omega}$ to $P$ by setting $\bigvee X\; :=\; sup (X)$ if $sup(X)\in P$. If $sup(X)\in S\setminus P$, then $\bigvee X$ is not defined. One calls the pair $(P,\bigvee)$ a {\it partial semilattice} \footnote{An {\it intrinsic} definition of $(P,\bigvee)$, i.e. one that avoids a comprising semilattice $S$ is possible. Essentially, each poset $(P,\le)$ can be enriched to a partial semilattice by focussing on any subsets $X_i\subseteq P$ for which $sup(X_i)$ happens to exist in $(P,\le)$, and by defining the partial operation just for these $X_i$. That may force the definition of additional $\bigvee Y$ but will not lead to contradictions.  Notice that \cite{gratzer} only considers {\it binary} partial operations. That would be too restrictive for us in the sequel.}. Let ${\cal R}$ be the set of semilattice  relations that bijectively correspond to those subsets $X\subseteq P$ for which $\bigvee X$ is defined. Then ${\cal R}$ is deductively closed since any relation not in ${\cal R}$ but deducible from ${\cal R}$, would in particular hold in $S$, but it doesn't.
Gathering all relations of type $a\ge b$ in ${\cal R}$  yields a unique poset $(P,\le)$. An order ideal $A$ of $(P,\le)$ is called {\it $\bigvee$-ideal} if 
\begin{displaymath}
\left(\forall X\subseteq A \right) \;\; \left(\bigvee X\in P \Rightarrow \bigvee X\in A\right)
\end{displaymath}
(Here of course "$\bigvee X\in P$" is shorthand for "the partial operation $\bigvee$ is defined for $X$".) One readily verifies that the set $Id(P,\bigvee)$ of all $\bigvee$-ideals is a closure system.

Let $(P,\bigvee)$ be a partial semilattice. The semilattice {\it freely generated} by $(P,\bigvee)$ is defined as $F_{\vee}\left(P,\bigvee \right) \; :=\; F_{\vee}\left(P,{\cal R}\right)$ where ${\cal R}$ is as above.

\begin{corol}\label{corol3}
Given a partial semilattice  $(P,\bigvee)$ the semilattice $F_{\vee}(P,\bigvee)$ is isomorphic to $Id(P,\bigvee)\setminus \{\emptyset\}$ with joins given as in (4). Moreover, $F_{\vee}(P,\bigvee)$ is up to isomorphism the unique semilattice $S$ such that 
\begin{itemize}
\item[(A)] There is a generating subset of $S$ which as a partial semilattice is isomorphic to $(P,\bigvee)$
\item[(B)] For each semilattice $T$ each partial morphism $\psi : P\ra T$ (i.e. $\psi(\bigvee X)=\bigvee\psi(X)$ whenever $\bigvee X$ exists) can be extended to a semilattice morphism $\Phi: S\ra T$.
\end{itemize}
\end{corol}

{\bf Proof:} By Theorem 1, $F_{\vee}(P,\bigvee)$ is isomorphic to $\mathscr{C}(\Sigma({\cal R}))\setminus\{\emptyset\}$  which clearly is $Id(P,\bigvee)\setminus\{\emptyset\}$. Properties $(A)$ and $(B)$ follow from $(i)$, $(ii)$ by taking into account that ${\cal R}$ is deductively closed (cf proof of Corollary 1). \hfill $\blacksquare$

\section{The core of a semilattice} \label{section5}
Each semilattice becomes a lattice by attaching a smallest element $0$. In this article we exploit the converse, namely that each lattice $L$ is in particular a semilattice, and that finite presentations  of semilattices are easier to handle. All of the sequel ows a lot to \cite{duquenne}. More details at the end of this section.
\newline

For starters, for any lattice $L$, the semilattice $L\setminus \{0\}$ is isomorphic to a finitely presented semilattice $F_{\vee}(P,{\cal R})$ when $P$ is taken as the set $J=J(L)$ of nonzero join irreducibles. As to ${\cal R}$, for $a\in L$ let $J(a):=\{p\in J\;|\; p\le a\}$, and for $X\subseteq J$ define its "natural closure" as $\ol{X} := J(\bigvee X)$.  If $\mathscr{C}$ is the associated closure system, then $a\mapsto J(a)$ yields an isomorphism from $L$ to $\mathscr{C}$. If $\Sigma$ is any implicational base of $\mathscr{C}$  and ${\cal R}$ the corresponding set of semilattice relations, then $L\setminus\{0\} \cong F_{\vee}(J,{\cal R})$  by Theorem 1.
\newline

Trouble is that $\Sigma$ is not readily found, and this requires some detours.
First, the set $J$ becomes a poset $(J,\le)$ with the ordering induced by $L$. With respect to the above closure system each closed set $X= J(a)$ is an order ideal, but (unless $L$ is distributive) not every order ideal $X$ of $(J,\le)$ is closed. 
 The trick will be to extend $J$ to a suitable partial semilattice $(P,\bigvee)$ and to compute $L$ as the closure system $Id(P,\bigvee)$ of all $\bigvee$-ideals of $(P,\bigvee)$ (see Corollary \ref{corol3}).
  \newline
  
Thus, given a lattice $L$, consider any subset $P\subseteq L$ for which
\begin{eqnarray}\label{eq11}
\psi : L\ra Id\left(P,\bigvee \right),\;\;\; a\mapsto P(a):=\{p\in P\;|\; p\le a\}
\end{eqnarray}
is a semilattice isomorphism. For instance, $P:=L$ does the job, but we strive for $P$ to be as small as possible. One always needs $J=J(L)\subseteq P$. Indeed, if we had $q\notin P$ for some $q\in J$, then $\psi(q)=\psi(q*)$ where $q*$ is the unique lower cover of $q$. This contradicts the injectivity of $\psi$. Conversely, for any superset $P$ of $J$, the map $\psi$ defined in (\ref{eq11}) will be injective because $a\ne b$ implies $J(a)\ne J(b)$ which, in view of $J\subseteq P$, amounts to $P(a)\cap J\ne P(b)\cap J$, and hence $P(a)\ne P(b).$   \newline

Assume again that $\psi$ in (\ref{eq11}) is an isomorphism. In order to derive a stronger necessary condition than just $J\subseteq P$, let us return to the quasi-closure of (6). For $q\in J$ one has $\{q\}^{\bullet} =\{q\}$ and hence $\{q\}^{\bullet}\ne \ol{\{q\}}=J(q)$ unless $q$ is an atom of $L$. In particular, each {\it non-atomic} $q\in J$ is such that $J(q)$ contains a nonclosed quasiclosed generating set.\newline

Consider now a {\it reducible} $a\in L$ such that $J(a)$ contains a quasiclosed nonclosed generating set $K$, that is, $K=K^{\bullet}\ne \ol{K}=J(a)$. We want to show that $a\in P$, and so by contraposition assume that $a\notin P$. 
By quasiclosedness  $K$ must be an order ideal of $(J,\le)$.  In order to see that 
\begin{eqnarray*}
K_P &:=& \left\{\bigvee X\;|\; X\subseteq K,\; \bigvee X\in P\right\}
\end{eqnarray*}
 is an order ideal of $(P,\le)$, fix $b\in P$ and $\bigvee X\in K_P$ with $b\le \bigvee X$. From $\bigvee X < \bigvee K$ (since $\bigvee K = a \notin P$) follows $J\left(\bigvee X\right)\subseteq K$ since $K=K^{\bullet}$. Hence there is $Y\subseteq J\left(\bigvee X\right)\subseteq K$ with $\bigvee Y = b$, and so $b\in K_P$. In fact, $K_P$ even is a $\bigvee$-ideal of the partial semilattice $(P,\bigvee)$ because if $\bigvee X$, $\bigvee Y \in K_P$ with $c:= \left(\bigvee X\right)\vee \left(\bigvee Y\right) \in P$, then $c= \bigvee\left(X\cup Y\right)\in K_P$ by definition of $K_P$. Hence $K_P$ is a member of $Id\left(P,\bigvee \right)$. But we claim it is not in the range of $\psi$. Since  $\bigvee K=a$, the only possibility for that to happen is $K_P=\psi(a)$.  However, by assumption there is some $r\in J(a)\setminus K$, and so $r\in \psi(a)\setminus K_P$. This contradiction shows that for $\psi$ to be an isomorphism, it is necessary that $P$ comprises the {\it join core}
 \begin{eqnarray} \label{eq12}
 K_{\vee}(L) &:=& J(L)\cup E_{\vee}(L),
 \end{eqnarray}  
 where the set of {\it join-essential} elements is
 \begin{eqnarray}\label{eq13}
 E_{\vee}(L) &:=&\{a\in L\;|\; \left( \exists K\subseteq J(a)\right) \; K^{\bullet} \ne \ol{K} = J(a)\}
 \end{eqnarray}
 
 \begin{thm}
 If $L$ is a finite lattice and $P :=K_{\vee}(L)$, then $L\cong F_{\vee}\left(P,\bigvee \right)$ as semilattices. Furthermore, for each subset $Q$ of $L$, one has $L\cong F_{\vee}(Q,\bigvee)$ if and only if $P\subseteq Q$.
 \end{thm}
 {\bf Proof:} By the deliberations above it suffices to prove the first claim. Recall from Corollary \ref{corol3} that $F_{\vee}(P,\bigvee)\simeq Id(P,\bigvee)\setminus \{\emptyset\}$. Because of $J\subseteq P$ we know that $\psi$ in (\ref{eq11}) is injective. In order to see that $\psi$ is surjective (and thus clearly an isomorphism), we show that for each nonempty $H$ in $Id\left(P,\bigvee\right)$ one has 
 \begin{eqnarray}\label{eq14}
  H\;=\;P(a), \;\;\;{\rm where}\;\;\; a\;\;:=\;\;\bigvee H \;\;\left(= \bigvee \left( J\cap H\right) \right).
 \end{eqnarray} 
  By contraposition, assume that $H$ is an inclusion-minimal counter example of (14). If we had $a := \bigvee H \in P$, then $a\in H$ because $H$ is a $\bigvee$-ideal, and so $P(a)\subseteq H$, i.e. $P(a)=H.$ Because by assumption $H\subsetneqq P(a)$, we conclude that $a\notin P$. Thus we get a desired contradiction to (\ref{eq12}), (\ref{eq13}) if we can establish $J\cap H$ as a quasiclosed nonclosed subset of $J(a)$. 
 \vskip 0.25cm
   As to "quasiclosed", we need to show that 
  $$ \left(X\subseteq J\cap H \;\;{\rm and}\; \bigvee X < a \right) \;\;{\rm implies}\;\; J\left(\bigvee X\right)\subseteq H.$$ 
  
 \ul{Case 1:}  $\bigvee X \in P$. Then $\bigvee X \in H$, and so $J\left(\bigvee X\right) \subseteq P\left(\bigvee X\right) \subseteq H$. \newline
 \ul{Case 2:} $\bigvee X \notin P$. Let $K_P$ be the $\bigvee$-ideal generated by $X$ (recall that $Id\left(P,\bigvee \right)$ is a closure system). Thus $X\subseteq K_P\subseteq H.$ In fact $K_P\subsetneqq H$ since $\bigvee X<a$, $\bigvee H=a$. But then (14) holds for $K_P$ by the minimality of $H$, and so $J\left(\bigvee X\right)\subseteq P\left(\bigvee X\right)=K\subseteq H$. 
 
 As to $J\cap H$ being nonclosed, let $b\in P(a)\setminus H$ be minimal. Suppose $b$ was reducible. Then $b=\bigvee Q$ with $Q\subseteq J(b)$, and $J(b)\subseteq H$ by the minimality of $b$. This is impossible since $Q\subseteq H$, $\bigvee Q\in P$ would force $\bigvee Q\in H.$ Therefore $b\in J(a)\setminus (J\cap H)$. \hfill $\blacksquare$

Duquenne \cite{duquenne} studies the {\it meet} core $K_{\wedge}(L)$ for various kinds of lattices. Up to duality his definition of $K_{\wedge}(L)$ matches our definition of the join core $K_{\vee}(L)$ in (\ref{eq12}). However, Duquenne's $\wedge$-essential elements by definition are the {\it reducible} members of $K_{\wedge}(L)$, whereas we defined $E_{\vee}(L)$ in such a way that $E_{\vee}(L)\cap J(L) =\{q\in J(L)\;|\; q \;{\rm is \;no \;atom}\}$. In so doing $E_{\vee}(L)$ matches the natural definition of  $E(\mathscr{C})$ in (a),(b) when $\mathscr{C}:= \{J(a)\;|\; a\in L\}$. For instance, for a finite Boolean lattice  $L$ equation (\ref{eq12}) becomes $K_{\vee}(L)=J(L)\cup \emptyset$, corresponding to the fact that only the {\it empty} family $\Sigma$ is a {\it nonredundant} implicational basis of $\mathscr{C}:=\{J(a)\;|\; a\in L\}\simeq 2^{J(L)}$. 

\section{The variety of all lattices}\label{section6}

In this and the next two sections we turn from freely generated semilattices to freely generated lattices.
Let $t_1, t_2, t_3$ be elements of any lattice. Then, clearly:
\begin{enumerate}
 \item [(i)]\;\; $t_1 \vee t_2 \;\leq \;t_3 \quad \LRa  \quad t_1 \;\leq \; t_3$ \;\; and \; \;$t_2 \;\leq \;t_3$
\item[(ii)] \;\; $t_3 \;\leq \; t_1 \wedge t_2 \quad \LRa \quad t_3 \;\leq t_1$ \;\; and \;\;$t_3 \;\leq \;t_2$
\item[(iii)] \;\; $t_3\; \leq \; t_1$ \;\; or \;\;$t_3 \;\leq \; t_2 \quad \Ra \quad t_3 \;\leq \; t_1 \vee t_2$
\item[(iv)] \;\; $t_1 \;\leq \;t_3$\;\; or \;\;$t_2 \;\leq \; t_3 \quad \Ra \quad t_1 \wedge t_2 \;\leq \; t_3$
\end{enumerate}
{\bf Definition:}\;
The {\it lattice} {\it freely generated by the poset} $P$ within a variety $\mathcal{V}$ of lattices is defined as  the up to isomorphism unique lattice $F\mathcal{V}(P,\le)$ in $\mathcal{V}$ which satisfies:
\begin{itemize}
\item[(a)]  $F\mathcal{V}(P,\le)$ is generated by $P$

\item[(b)]   For each lattice $T\in \mathcal{V}$ and each monotone map $\phi : P \ra T$, there is a homomorphism $\Phi : F\mathcal{V}(P,\le) \ra T$ that extends $\phi$.
\end{itemize}

\begin{figure}[!h]
\begin{center}
\psfrag{fr}{$FL(P,\le)$\;\; = }
 \includegraphics[scale=.4]{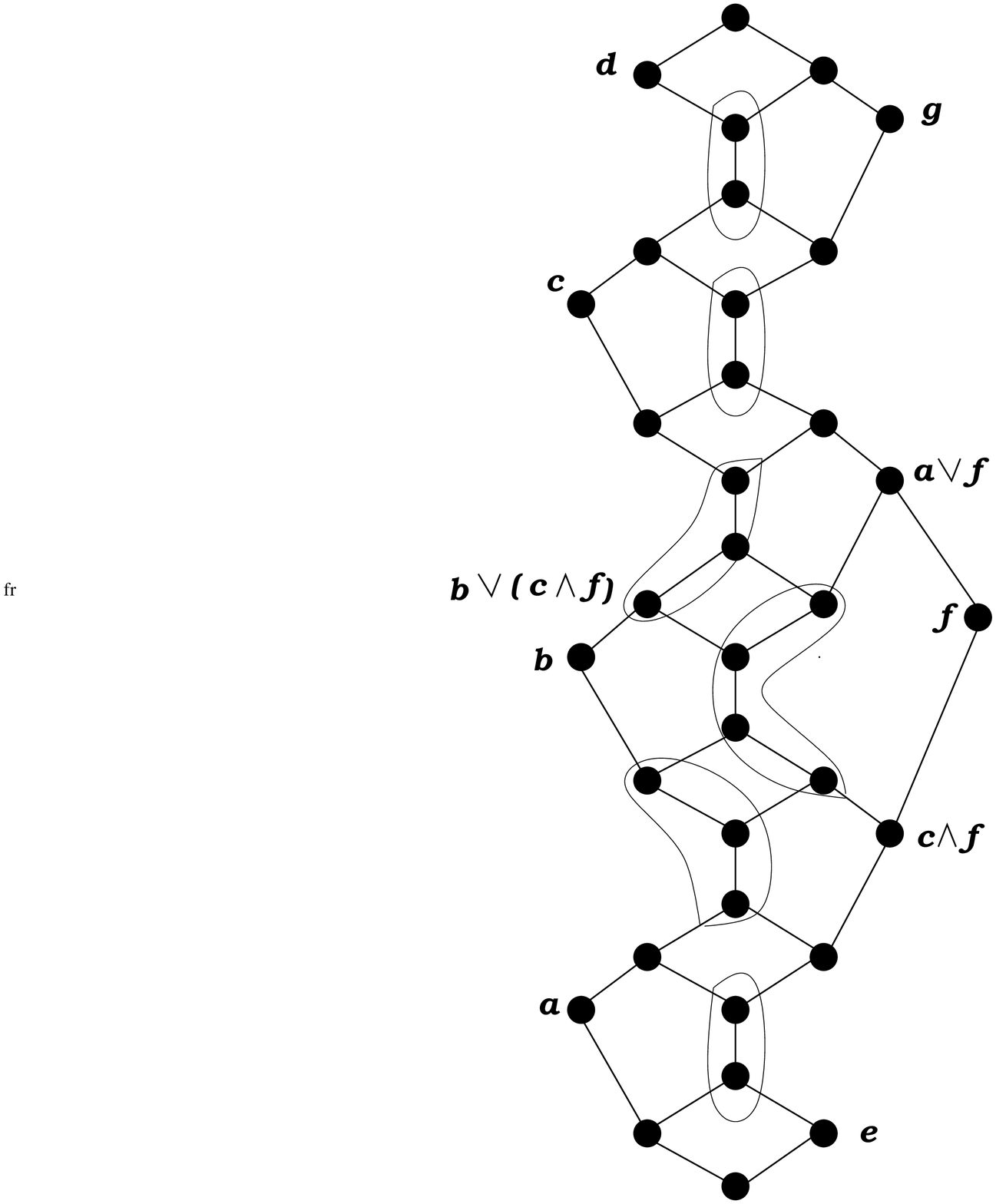}
\caption{}
\label{figg4}
\end{center}
\end{figure}

In this section $\mathcal{V}$ is the variety of {\it all} lattices and we write $FL(P,\le)$ for $F \mathcal{V}(P,\le)$.
 It turns out\footnote{For unordered sets $P$ this has been proven by P. M. Whitman in 1941, although J. B. Nation 
  states that Skolem achieved essentially the same in 1920. The natural extension to posets $(P,\le)$ is due to R. P. Dilworth 1945. Reproving a 1964 result of R. A. Dean, H. Lakser \cite{lakser} dealt with the lattice freely generated by $(P,\le)$ and additionally preserving certain finite joins and meets in $(P,\le)$.} that the converse implication in (iii) respectively (iv) holds in every lattice $FL(P,\le)$. 
 \vskip 0.25cm
 For instance, $P$ being the poset from figure \ref{figure3}$(a)$, is it true that $b \vee (c \wedge f) \leq c \wedge (a \vee f)$ in $FL(P,\le)$? By (i) the truth amounts to
$$b \leq c \wedge (a \vee f) \quad \mbox{and} \quad (c \wedge f ) \leq c \wedge (a \vee f),$$
which by (ii) amounts to
$$b \leq c \quad \mbox{and} \quad b \leq a \vee f \quad \mbox{and} \quad c \wedge f \leq c \quad \mbox{and} \quad c \wedge f \leq a \vee f.$$
Now by the converse of (iii), $b \not\leq a \vee f$ since $b \not\leq a, b \not\leq f$ in $P$. Therefore also $b \vee (c \wedge f) \not\leq c \wedge (a \vee f)$. For our $(P,\le)$ of figure \ref{figure3}$(a)$, $FL(P,\le)$ happens to be finite and is depicted in figure \ref{figg4}. Observe that the "bubbles" in figure \ref{figg4} are the congruence classes of the epimorphism $FL(P,\le)\ra FD(P,\le)$ where $FD(P,\le)$ is the free distributive lattice discussed in section \ref{section7}. The lattices $FL(P,\le)$ will recur in section \ref{section6} of Part II.

\section{The variety of distributive lattices}\label{section7}

If ${\cal V}$ is the variety ${\cal D}$ of distributive lattices, then $F{\cal V}(P,\le)$ will be written as $FD(P,\le)$. Like every distributive lattice, $FD(P,\le)$  is a subdirect product of two element factor lattices $D_2=\{0,1\}$. These factors can be neatly distinguished by the set of elements ${\cal F}\subseteq P$ that map upon $1$ (as opposed to $0$) under the projection $\pi : FD(P,\le)\ra D_2$. Since $\pi$ is monotone, ${\cal F} = \pi^{-1}(1)$ is a (nonempty) order filter of $(P,\le)$. {\it All} order filters ${\cal F}$ arise in this way by the universal mapping property of $FD(P,\le)$.  Thus we find that $FD(P,\le)$ will have these $12$ subdirectly irreducible factors:  

\newpage

\begin{figure}[!h]
\begin{center}
\psfrag{a}{$a$} \psfrag{b}{$b$} \psfrag{c}{$c$} \psfrag{d}{$d$} \psfrag{e}{$e$} \psfrag{f}{$f$} \psfrag{g}{$g$} 
\includegraphics[scale= .4]{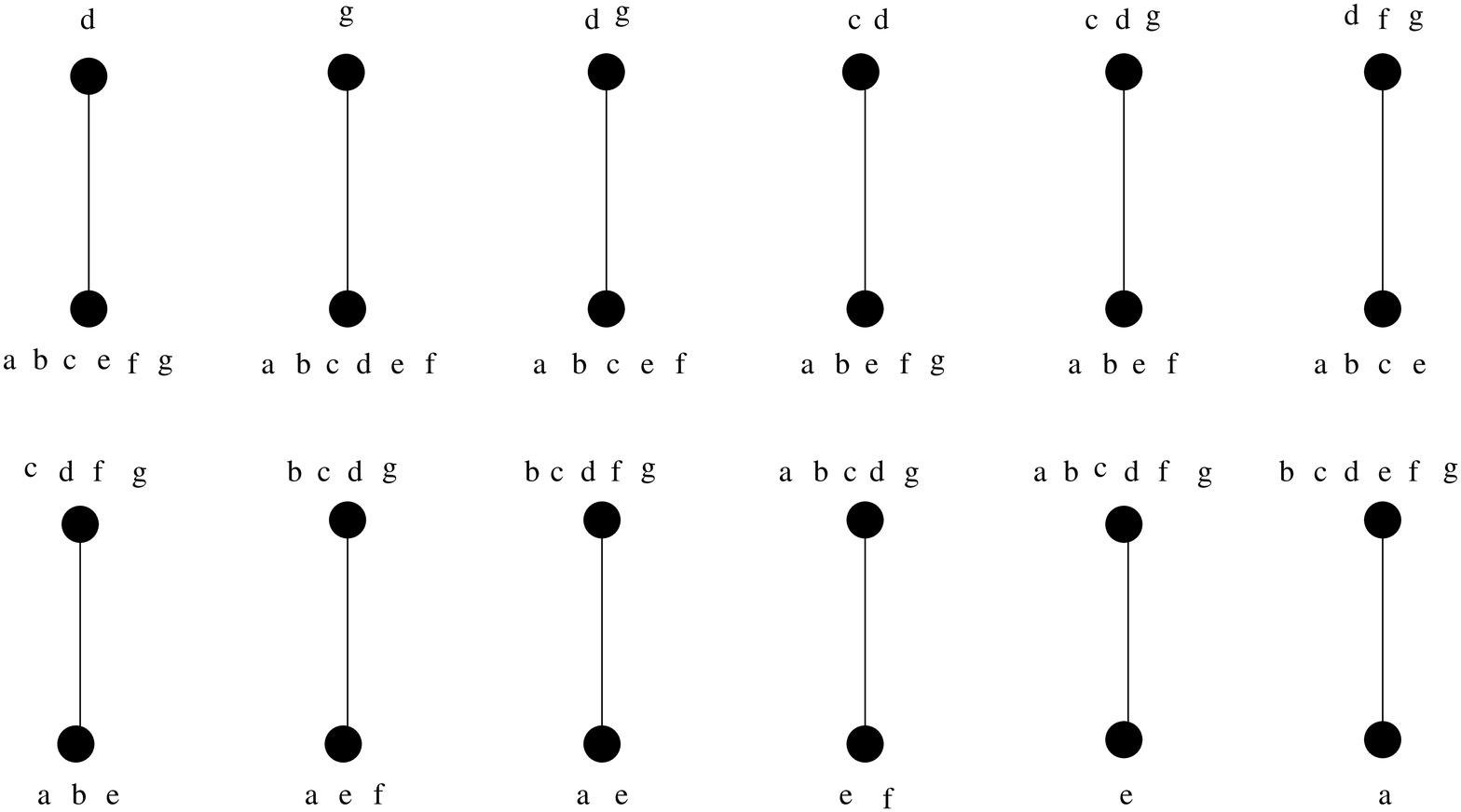}
\caption{}
\label{figur5}
\end{center}
\end{figure}

 At this point we are stuck because the information provided by the "$P$-labellings" in Fig.\ref{figur5} is not enough to construct $FD(P,\le)$. What one needs are the {\it connection maps} between any two $P$-labellings.  We shall persue this line of thought in our sequel paper which deals with arbitrary finitely generated varieties of lattices (of which ${\cal D}$ is the simplest example). But here we tackle $FD(P,\le)$ in a way that avoids subdirect products. The core is contained in the next lemma.
\begin{lm}\label{lemma1}
 Let $F$ and $T$ be distributive lattices such that $F$ is finite and $J_0(F) := J(F)\cup\{0\}$ is closed under meets. Then each $\wedge$-preserving map $\overline{\phi}: J_0(F)\rightarrow T$ can be extended to a homomorphism $\Phi: F \ra T$.
\end{lm}
{\bf Proof:} For all $x\in F$ put 
$$\Phi(x) \quad := \quad \bigvee \{\overline{\phi}(s)| s\le x\}$$

where by convention $s$ (or $t$ or $r$) ranges over $J_0(F)$. Because $\overline{\phi}$ is monotone, $\Phi$ extends $\overline{\phi}$ and is monotone itself. It thus remains to show that $\Phi$ preserves $\wedge$ and $\vee$. As to $\wedge$, by assumption $s,t\in J_0(F)$ implies  $s\wedge t\in J_0(F)$, and $\overline{\phi}(s)\wedge \overline{\phi}(t)=\overline{\phi}(s\wedge t)$. Hence 

\begin{eqnarray*}
 \Phi(x)\wedge \Phi(y) & = & \bigvee \{ \overline{\phi}(s)| s\le x\} \wedge \bigvee \{\overline{\phi}(t)| t\le y\}  \\
 & = & \bigvee \{\overline{\phi}(s) \wedge \overline{\phi}(t)| s\le x, t\le y\} \qquad {\rm by\; distributivity}\\
 & = & \bigvee \{\overline{\phi}(s\wedge t)| s\le x, t\le y\} \\
 & = & \bigvee \{\overline{\phi}(r)| r\le x\wedge y \} \;\; =\;\; \Phi(x\wedge y)\\
\end{eqnarray*}

By distributivity each $s\in J_0(F)$ is join-prime, that is, $s\le x\vee y$ implies $s\le x$ or $s\le y$. Hence 

\begin{eqnarray*}
 \Phi(x\vee y) &=& \bigvee \{\overline{\phi}(s)| s\le x\vee y\}\\
                &=& \bigvee \{\overline{\phi}(s)| s\le x \;\;{\rm or}\;\; s\le y\}\\
               &\le& \bigvee \{\overline{\phi}(s)| s\le x\} \vee \bigvee \{\overline{\phi}(s)| s\le y\}\\
               &=& \Phi(x) \vee \Phi(y)
\end{eqnarray*}

The inequality $\ge$ is trivial. \hfill $\blacksquare$

 \begin{thm} \label{thme3}
 The lattice $FD(P,\le)$ freely generated by the finite poset $(P,\le)$ within the variety
  of all distributive lattices  is isomorphic to $F_{\vee}\left(F_{\wedge}(P,\le),\le \right)$.
 \end{thm}
 
{\bf Proof:} For elements $p,q$ of any poset, if the meet $p\wedge q$ happens to exist, then one checks that
\begin{eqnarray*}
(p\wedge q)\!\downarrow &=& p\!\downarrow \cap \;q\!\downarrow.
\end{eqnarray*}
We view the $\wedge$-semilattice $J_0 := F_{\wedge}(P,\le)$ as a poset with smallest element $0$. By Corollary 1 the $\cup$-semilattice of nonempty order ideals of $J_0$ can be identified with $F:=F_{\vee}(J_0,\le)$. Since   $F$ has a smallest element (corresponding to $\{0\}$), it is a distributive lattice with $J_0(F)$ equal to $J_0$. By the above remark $J_0$ is closed under meets.

Let now $T$ be any distributive lattice and $\phi: P\ra T$ a monotone map. Since $J_0$ is the free $\wedge$-semilattice generated by $(P,\le)$, $\phi$ can be extended to a $\wedge$-preserving map $\overline{\phi}$ on $J_0$. By Lemma \ref{lemma1}, $\overline{\phi}$  further extends to a homomorphism $\Phi : F\ra T$. \hfill $\blacksquare$

The proof given here is believed to be new. Theorem \ref{thme3} is a special case of  more general (but  clumsier) results, e.g. by W. R. Tunnicliffe \cite{tunnicliffe} or Yongming Li \cite{yongming}. 

\begin{ex}{\rm Consider the poset $(P,\le)$ from Fig.\ref{figure3}($a$). One checks that the fat subset $J_0$ of join irreducibles (including $0$) of $FD(P,\le)$ is a meet subsemilattice, and it is isomorphic to $F_{\wedge}(P,\le)$ from Fig.\ref{figure3}($c$). The elements $a,\cdots, g$ of $P$ correspond to the doubly-irreducible elements of $FD(P,\le)$:}\end{ex}

\newpage
\begin{figure}[!h]
\begin{center}
\psfrag{t}{$FD(P,\le)\; \;= $}
 \includegraphics[scale=.4]{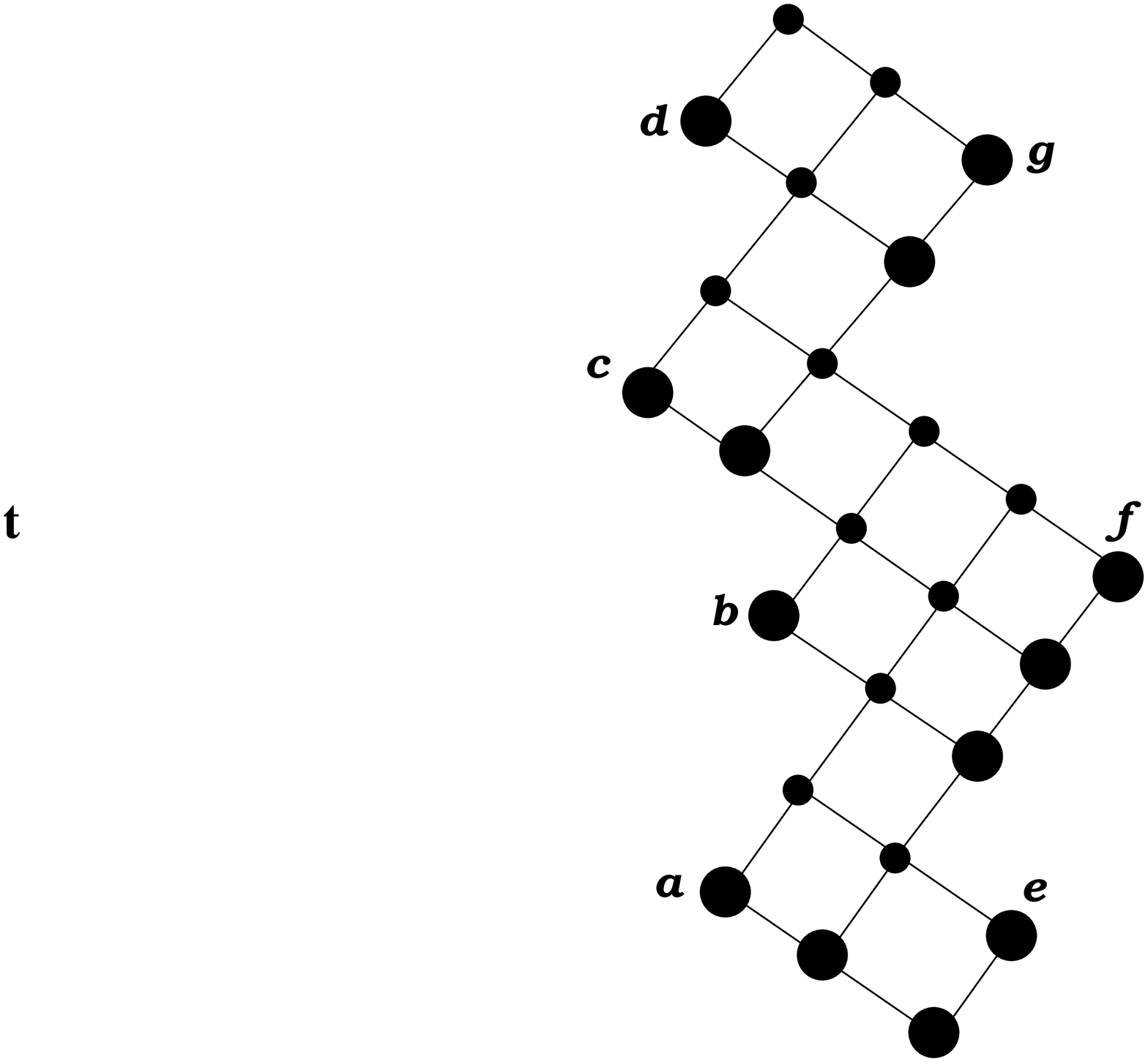}
\end{center}
\caption{}
\end{figure}
The subdirectly irreducible factors $\pi : D \ra D_2$ of every finite distributive lattice $D$ bijectively correspond to the join irreducibles $q$ in that 
\begin{eqnarray}
\pi^{-1}(1) &=& \{x\in D\;|\;x\ge q\} \;\;\; {\rm and}\;\;\; \pi^{-1}(0) \;\;=\;\;\{x\in D\;|\;x\ngeq q\}.
\end{eqnarray}
If for $D=FD(P,\le)$ we take (say) $q=b$, we have $\pi^{-1}(1)\cap P=\{b,c,d,g\}$ and $\pi^{-1}(0)\cap P=\{a,e,f\}$  corresponds indeed to the eigth $P$-labelling in Fig.\ref{figur5}.

{\it Open problem:} Set up an algorithm that computes {\it arbitrary} finitely presented lattices within $\cal{D}$, such as the free distributive lattice generated by $a,b,c,d$ and subject to say $a\vee b=(c\wedge d)\vee (a\wedge b)$.

\section{The variety of Boolean algebras}\label{section8}

The class of all Boolean lattices is no variety since it is not closed under taking sublattices. But it becomes  a variety if $0,1$ are elevated to nullary operations (constants) and complementation $x \mapsto {\ol x}$ is added as unary operation. One then speaks of Boolean {\it algebras} and  accordingly the definition in section \ref{section6} has to be adjusted to the extent that for the {\it free Boolean algebra $FB(P,\le)$ generated by the poset} $(P,\le)$, the homomorphism $\Phi$ in $(b)$ must be a {\it Boolean} homomorphism in the sense that $\Phi (0) =0, \Phi (1) =1$, and $\Phi ({\ol x}) = {\ol {\Phi(x)}}$ for all $x \in FB(P,\le)$.  This section is based on [6, p. 107]. 
 
In order to describe $FB (P,\le)$ for finite $P$,  observe that any Boolean algebra $B(x_1, \cdots,$ $x_s)$  which is generated by $s$ elements of some comprising Boolean algebra, by distributivity and De Morgan's laws equals
\begin{eqnarray}\label{eq16}
B(x_1, \cdots, x_s) \quad = \quad \left\{ \left. \ds\bigvee_{\delta \in I} \left(x^{\delta_1}_1 \wedge x^{\delta_2}_2 \wedge \cdots \wedge x^{\delta_s}_s \right)\right| \ I \subseteq \{0,1\}^s \right\}.
\end{eqnarray}
Here $\delta = (\delta_1, \delta_2, \cdots, \delta_s)$ is a $0,1$-vector and $x^1_i := x_i, \ x^0_i : = \ol{x_i}$. From (\ref{eq16}) we conclude that the $t$ atoms of $B(x_1, \cdots, x_s)$ are exactly the {\it nonzero} elements among the at most $2^s$ elements $x^{\delta_1}_1 \wedge \cdots \wedge x^{\delta_s}_s$. It follows that if $FB(s)$ is the free Boolean algebra generated by $s$ unordered elements, then
\begin{eqnarray}\label{eq14}
\left| FB(s)\right| \quad \le \quad 2^{(2^s)}.
 \end{eqnarray}

\begin{ex} \label{example7}
{\rm Let $X:= \{a, b, \cdots, i\}$. What is the number $t$ of atoms of the Boolean algebra $B (A_1, A_2, A_3, A_4) \subseteq {\cal P}(X)$ generated by $A_1 := \{a, c, d, f\}, \ A_2 := \{b, e, f, i\}, \ A_3 := \{a, b, c, d, e, g, i\}$ and $A_4 := \{e, g, h\}$?}
\end{ex}

\hspace*{3cm} \begin{tabular}{c|c|c|c|c|c|c|c|c|c|} 
& $a$ & $b$ & $c$ & $d$ & $e$ & $f$ & $g$ & $h$ & $i$\\ \hline
$A_1 =$ & $1$ & $0$ & $1$ & $1$ & $0$ & $1$ & $0$ & $0$ & $0$ \\ \hline 
$A_2=$ & $0$ & $1$ &  $0$ & $0$ & $1$ & $1$ & $0$ & $0$ & $1$ \\ \hline
$A_3=$ & $1$ & $1$ & $1$ & $1$ & $1$ &$0$ & $1$ &$0$ & $1$\\ \hline
$A_4=$ & $0$ & $0$ & $0$ & $0$ & $1$ & $0$ & $1$ & $1$ & $0$ \\ \hline \end{tabular}

Here the characteristic vectors of the subsets $A_i \subseteq S$ are listed as the rows of a $4 \times 9$ matrix. Then the number $t$ of equivalence classes of equal columns obviously is the number of $(\delta_1, \delta_2,\delta_3, \delta_4)$ with $A^{\delta_1}_1 \cap A^{\delta_2}_2 \cap A^{\delta_3}_3\cap A^{\delta_4}_4 \neq \emptyset$. For instance, the three equal columns labelled by $a, c, d$ yield $(\delta_1, \delta_2, \delta_3, \delta_4) = (1, 0, 1, 0)$ and $A^1_1 \cap A^0_2 \cap A_3^1 \cap A^0_4 = \{a, c, d\} \neq \emptyset$. In our case $t=6$ whence $|B (A_1, A_2, A_3, A_4)| = 2^6 = 64$.  \hfill $\square$

Picking $s$ random sets from a $r$-set amounts to pick $r$ random $0,1$-columns of length $s$. In Example \ref{example7} we had $r=9$ and $s=4$. As we saw, had the columns be distinct (which doesn't imply distinct rows), then $A_1,\cdots, A_s$ would have generated ${\cal P}(X)$. The probability of getting distinct columns is easily calculated:

\begin{thm}
 Let $X$ be a $r$-set. Pick $s$ (not necessarily distinct) sets $A_i \in {\cal P}(X)$ at\\ random. The probability that $B(A_1, \cdots, A_s)$ equals ${\cal P}(X)$ is 
 $$\ds\frac{2^s-1}{2^s} \cdots \frac{2^s-2}{2^s} \cdots \frac{2^s - (r-1)}{2^s}$$
\end{thm}
For instance, the probability that $s=5$ random subsets $A_i \subseteq X := \{1, 2, 3, 4, 5\}$ generate the whole powerset is $0.72$. When all $A_i \subseteq \{1, 2, \cdots, 10\}$, the probability is $0.21$. When $r=2^s$, the probability for $| B (A_1, \cdots, A_s)| = |{\cal P}(X)| = 2^{(2^s)}$ is still $> 0$.  In fact, in view of (\ref{eq14}) such a $B(A_1,\cdots, A_s)$ must be isomorphic to $FB(s)$. Let us now turn from $s$-element antichains to general posets $(P,\le)$.

\begin{thm}
If $(P,\le)$ is finite, then $FB (P,\le)$  has exactly $t$ atoms, where $t$ is the
 number of 
order filters $\emptyset \subseteq F \subseteq P$. 
\end{thm}

{\bf Proof:} Assume w.l.o.g. $(P, \leq)$ is $(\{1, 2, \ldots, s\}, \leq )$ and that $(\{A_1, \cdots, A_s\}, \subseteq )$ is any fixed set system such that $i \leq j \Leftrightarrow A_i \subseteq A_j$ for all $i, j \in P$. In view of Example \ref{example7} observe that $A^{\delta_1}_1 \cap \cdots \cap A^{\delta_s}_s$ can be nonempty only if ``$\delta_i = 1 \Ra \delta_j =1$'' whenever $i \leq j$. Hence the number of vectors $\delta = (\delta_1, \cdots, \delta_s)$ with $A^{\delta_1}_1 \cap \cdots \cap A^{\delta_s}_s \neq \emptyset$ is {\it at most} the number $t$ of order filters of $P$. On the other hand, choosing $(\{A_1, \cdots, A_s \}, \subseteq )$ appropriately (as in Example 8 below), one can indeed obtain $t$ nonempty sets $A^{\delta_1}_1 \cap \cdots\cap A^{\delta_s}_s$.  This Boolean algebra of largest cardinality generated by $(P,\le)$ can  only be $FB(P,\le)$ \hfill $\blacksquare$

\begin{ex} {\rm What is the size of the free Boolean algebra $FB (P,\le)$ generated by our companion poset $(P, \leq)$ of figure \ref{figure3}$(a)$?
As opposed to Example \ref{example7}, here the $0,1$-matrix is generated by concatenating the {\it columns}. Namely, if we list the characteristic vectors of the $t=14$ order filters of $(P, \leq)$ as the (distinct) columns of a $7 \times 14$ matrix, then the sets corresponding to the rows yield a set system isomorphic to $(P, \leq)$ (check), and all intersections $A^{\delta_a}_a \cap \cdots \cap A^ {\delta_g}_g$, where $\delta$ is a column of the matrix, are nonempty (as argued in Example 7). Thus $|F B (P,\le)| = 2^{14} = 16384$.}
\end{ex}

\begin{tabular}{c|c|c|c|c|c|c|c|c|c|c|c|c|c|c|} 
& $\emptyset$ & $a\!\uparrow$ & $b\!\uparrow$ & $c\!\uparrow$ & $d\!\uparrow$ & $e\!\uparrow$ & $f\!\uparrow$ & $g\!\uparrow$ & $a,e\! \uparrow$ & $a,f\!\uparrow$ & $b,f \!\uparrow$ & $c,f\!\uparrow$ & $c,g\!\uparrow$ & $d,g\!\uparrow$\\ \hline
$A_a=$ & $0$ & $1$ & $0$ & $0$ & $0$ & $0$ & $0$ & $0$  & $1$ & $1$ & $0$ & $0$ & $0$ & $0$ \\ \hline
$A_b=$ & $0$ & $1$ & $1$ & $0$ & $0$ & $1$ & $0$ & $0$ & $1$ & $1$ & $1$ & $0$ & $0$ & $0$ \\ \hline
$A_c=$ & $0$ & $1$ & $1$ & $1$ & $0$ & $1$ & $0$ & $0$ & $1$ & $1$ & $1$ & $1$ & $1$ & $0$ \\ \hline
$A_d=$ & $0$ & $1$ & $1$& $1$ & $1$ & $1$ & $1$ & $0$ & $1$ & $1$ & $1$ & $1$ & $1$ & $1$ \\ \hline
$A_e =$ & $0$ & $0$ & $0$ & $0$ & $0$ & $1$ & $0$ & $0$ & $1$ & $0$ & $0$ & $0$ & $0$ & $0$ \\ \hline
$A_f=$ & $0$ & $0$ & $0$ & $0$ & $0$ & $1$ & $1$ & $0$ & $1$ & $1$ & $1$ & $1$ & $0$ & $0$ \\ \hline
$A_g =$ & $0$ & $1$ & $1$ & $0$ & $0$ & $1$ & $1$ & $1$ & $1$ & $1$ & $1$ & $1$ & $1$ & $1$ \\ \hline \end{tabular}

\end{document}